\documentclass[12pt]{article}
 
\usepackage[margin=1in]{geometry} 
\usepackage{amsmath,amsthm,amssymb,graphicx,mathtools,tikz,hyperref}
\usetikzlibrary{positioning}
\usepackage{hyperref}
\usepackage[bottom]{footmisc}
\usepackage[final]{pdfpages}

\newcommand{\footremember}[2]{%
    \footnote{#2}
    \newcounter{#1}
    \setcounter{#1}{\value{footnote}}%
}
\newcommand{\footrecall}[1]{%
    \footnotemark[\value{#1}]%
}

 \hypersetup{
 colorlinks,
 linkcolor=blue
 }
\begin{document}
\date{\today}

\title{Analytic Formulas for the Sample Moments of the Sample Correlation over Permutations of Data}


\author{
  Marc Jaffrey\footremember{uw}{University of Washington (\href{mailto:mjaffrey@uw.edu}{mjaffrey@uw.edu})}\footremember{rf}{RootFault LLC} \and
  Michael Dushkoff\footremember{rit}{Rochester Institute of Technology (\href{mailto:mad1841@rit.edu}{mad1841@rit.edu})}\footrecall{rf}
}
 
\maketitle

\section*{Introduction}
Presented is a branching inductive formula for deriving analytic formulas for the sample moments of the distribution of the sample correlation coefficient $r$, 
\begin{equation}
    r = \frac{\sum_i (x_i -\mu_x)(y_i-\mu_y)}{ \sqrt{\sum_k(x_k-\mu_x)^2\sum_j(y_j-\mu_y)^2}},
    \label{ssc}
\end{equation}
over permutations of the data in terms of the central moments of the data, denoted
\begin{equation}
    \langle r_{\Pi}^k\rangle = \frac{1}{|\Pi|} \sum_{\pi \in \Pi} r_{\pi}^{k}
    \label{mdef}
\end{equation}
for $k \in \mathbb{N}$, where $\Pi$ is the set of all possible permutations on the data and $r_{\pi}$ denotes the sample correlation over the permutation $\pi \in \Pi$ of the data.

While these exact formulas are interesting by themselves, placed in a proper statistical framework they could open the door to the possibility of more precise and computationally efficient methods of evaluating the p-value for a hypothesis test of Pearson's correlation \cite{pearson1931test,pearson1895correlation}. %

%
\section*{Inductive Formula for $\langle r_{\Pi}^k\rangle$}
\subsection*{Notation} Given a dataset $D_n = \lbrace ( x_i,y_i )|\ i=1,..,n \rbrace \subset \mathbb{R}^2$, let $\Pi=Perm(n)$ be the set of all permutations on $\{ 1,\dots,n\}$. For $\pi=(\pi_1,\dots,\pi_n)\in \Pi$, defined such that $i\underset{\pi}{\to} \pi_i$, define
\begin{equation}
 \pi(D_n)=\lbrace(x_i,y_{\pi_i})|i=1,..,n\rbrace, 
 \label{py}
\end{equation}
where the permutation $\pi$ acts on the $y$-coordinate alone. Denote Pearson's sample correlation over a permutation $\pi$ of the data as $r_\pi = r(\pi(D_n))$
and define 
\begin{eqnarray}
   r_{\Pi}: [-1,1] \to [0,1],
\end{eqnarray}
the distribution of $r$ over the permutations of $D_n$.
Lastly, let $\hat{x}_i \equiv (x_i -\mu_x)$ and $\hat{y}_i \equiv   (y_i - \mu_y)$.

\subsection*{Main Result}
Given $D_n$, fix $k$. The $k^{th}$ moment of $r_{\Pi},$  is given by
\begin{equation}
    \langle r_{\Pi}^k \rangle   =   \frac{1}{n^k n!\hat{\sigma}_x^k \hat{\sigma}_y^k} \sum_{m=1}^{k} \ \ \sum_{n_1+...+n_m=k} \binom{k}{n_1,..,n_m}^* X_{(n_1,..,n_m)}^{m,k} (n-m)!Y_{(n_1,..,n_m)}^{m,k}\cdot h_{n,m}
    \label{pcpm}
\end{equation}
$X_{(n_1,..,n_m)}^{m,k} \text{ and } Y_{(n_1,..,n_m)}^{m,k} $ are branching inductive formulas, such that for $1<m\leq k$ and for fixed non-zero positive integers $n_1+...+n_m=k$ we have:
\begin{equation}
    X_{(n_1,..,n_m)}^{m,k}  =   \ n \langle  \hat{x}^{n_m} \rangle X_{(n_1,...,n_{m-1})}^{m-1,k-n_m}  - \sum_j X_{(n_1 + n_m\delta_{ij},...,n_{m-1}+n_m\delta_{(m-1)j})}^{m-1,k} 
    \label{inducx}
\end{equation}
where $\delta_{ij}$ is the standard delta function, and similarly defined
\begin{equation}
    Y_{(n_1,..,n_m)}^{m,k}  =   \ n \langle  \hat{y}^{n_m} \rangle Y_{(n_1,...,n_{m-1})}^{m-1,k-n_m}  - \sum_j Y_{(n_1 + n_m\delta_{ij},...,n_{m-1}+n_m\delta_{(m-1)j})}^{m-1,k}
    \label{inducy}
\end{equation}
For $m=1$, we have
\begin{equation}
X_{(k)}^{1,k}  =  n\langle \hat{x}^k \rangle  \text{ and }
     Y_{(k)}^{1,k} =  n\langle \hat{y}^k \rangle
\end{equation}
The stared multinomial coefficient
\begin{equation}
    \binom{k}{n_1,\dots,n_m}^* = 
\frac{1}{d_1!\cdots d_r!}\binom{k}{n_1,\dots,n_m}
\label{multi}
\end{equation}
is an adjustment of the usual multinomial coefficient accounting for degeneracy in $n_1,...,n_m$. Partitioning $n_1,...m_m$ into subsets, $g_1,..,g_r$, by the equivalence relation $n_i\equiv n_j \iff n_i=n_j$, then $d_i=|g_i|$. This degeneracy leads to an over counting represented by the multinomial coefficient which is correct by dividing out $d_1!\cdots d_r!$.

The term $h_{n,m}$ accounts for the inability to compute higher order terms in the sum when the number of data points is less than the moment order being computed,  by setting them to zero in the formula,
\begin{equation}
h_{n,m} = \left\{ \begin{array}{rcl}
 0 & \mbox{for}
& n-m < 0 \\ 1 & \mbox{for} & n-m \geq 0
\end{array}\right. 
\end{equation}
Lastly, for notational simplicity, in \eqref{pcpm} the following convention is employed,
\begin{equation}
    \hat{\sigma}_z = \sqrt{\frac{1}{n}\sum (z_i -\mu_z)^2}
    \label{sigma}
\end{equation}

\section*{Exact Formulas for $k=1,..,5$}
\label{app:form}
\noindent
From the induction formula every moment can be analytically determined. Below are the exact formulas for the first five moments in terms of the easily computed moments of the data. As in the main body of the paper, $\hat{\sigma}_x$ and $\hat{\sigma}_y$ are defined by \eqref{sigma}. Additionally, we employ the following notation for simplicity sake with $\chi_k=\langle x^k\rangle$, $\nu_k=\langle y^k \rangle$, and $\mu_{k,j}=\langle x^k \rangle \langle y^j \rangle$, all central moments. The first five moments of the sample correlation are as follows:
\begin{eqnarray*}
    \langle r_{\Pi}^1 \rangle  = &\  0
    \\
    \langle r_{\Pi}^2 \rangle = &\ \frac{1}{(n-1)}\\
    \\
     \langle r_{\Pi}^3 \rangle  = &\  \frac{\mu_{3,3}}{\hat{\sigma}_x^3\hat{\sigma}_y^3}\bigg[\frac{1}{n^2}h_{n,1} +\frac{3}{n^2(n-1)}h_{n,2} +\frac{4}{n^2(n-1)(n-2)}h_{n,3}\bigg]
\\
\\
      \langle r_{\Pi}^4 \rangle  = & \ \frac{1}{\hat{\sigma}_x^4 \hat{\sigma}_y^4} \bigg[\frac{\mu_{4,4}}{n^3}h_{n,1} + 
        \frac{4\chi_4 \nu_4}{n^3(n-1)}h_{n,2} +
          \frac{3[n^2\hat{\sigma}_x^4-n\chi_4][n^2\hat{\sigma}_y^4 -n\nu_4]}{n^5(n-1)}h_{n,2}
        \\ \\
        & +\frac{6[2n\chi_4 -n^2\hat{\sigma}_x^4][2n\nu_4 -n^2\hat{\sigma}_y^4]}{n^5(n-1)(n-2)}h_{n,3}
        + \frac{9[2n\chi_4 - n^2 \hat{\sigma}_x^4][2n\nu_4 - n^2 \hat{\sigma}_y^4]}{n^5(n-1)(n-2)(n-3)}h_{n,4}\bigg] 
\end{eqnarray*}
\begin{eqnarray*} 
       \langle r_{\Pi}^5 \rangle  = & \ \frac{1}{\hat{\sigma}_x^5 \hat{\sigma}_y^5} \bigg[ \frac{\mu_{5,5}}{n^4}h_{n,1} 
         + 5\frac{\mu_{5,5}}{n^4(n-1)}h_{n,2}
         +
         10\frac{[n^2\chi_3 \chi_2 -n\chi_5][n^2\nu_3 \nu_2 -n\nu_5]}
         {n^6(n-1)}h_{n,2}
         \\ 
         \\ 
         &
         +10\frac{[2n\chi_5 -n^2\chi_3\chi_2][2n\nu_5 -n^2\nu_3 \nu_2]}{n^6(n-1)(n-2)}h_{n,3}
         \\
         \\
         & + 60\frac{[n\chi_5-n^2\chi_3\chi_2][n\nu_5-n^2 \nu_3 \nu_2]}{n^6(n-1)(n-2)} h_{n,3}\\
         \\
         & +10\frac{[6n\chi_5 - 5n^2\chi_3\chi_2][6n\nu_5 - 5n^2\nu_3\nu_2]}{n^6(n-1)(n-2)(n-3)}h_{n,4
         }\\
         \\& +
         \frac{16[6n\chi_5 - 5n^2\chi_3\chi_2][6n\nu_5 - 5n^2\nu_3\nu_2]}{n^6(n-1)(n-2)(n-3)(n-4)}h_{n,5}\bigg]
\end{eqnarray*}

\subsection*{Validation}
\label{app:valid}
In order to demonstrate the validity of these formulas, we randomly generated datasets of sizes $n=\{3,\dots,8\}$ over $100$ trials at each $n$. We compared the derived moments to the empirical moments computed directly from $\Pi(D_n)$, which can be done in the case of small datasets, however for larger datasets ($n>8$) this becomes computationally impractical.
\\
\\
The mean squared error was computed as:
\begin{equation}
  MSE = \frac{1}{N_{trials}}\sum \left(\langle r^k\rangle_{\Pi(D_n)} - \langle r^k\rangle_{exact}\right)^{2}
\end{equation}
\\
This validation procedure was performed using MATLAB with double-precision floating point computations. The errors tabulated in Table \ref{tab:validation} are within machine epsilon error indicating that the formulas are indeed exact.

\begin{table}[h]
\centering
\begin{tabular}{| c || c | c | c | c | }

\hline
 &\multicolumn{4}{|c|}{MSE of $k^{th}$ Moment } \\ [0.5ex] 
 
  Sample Size & 2&3&4&5 \\ [0.5ex] 
 
 \hline
 $n=3$ & $1.20\mathrm{e}{-32}$ & $2.57\mathrm{e}{-32}$ & $3.99\mathrm{e}{-32}$ & $4.82\mathrm{e}{-32}$\\ 
 \hline
 $n=4$ & $6.81\mathrm{e}{-33}$ & $2.46\mathrm{e}{-33}$ & $1.11\mathrm{e}{-32}$ & $3.18\mathrm{e}{-33}$\\ 
 \hline
 $n=5$ & $6.06\mathrm{e}{-33}$ & $1.23\mathrm{e}{-33}$ & $5.75\mathrm{e}{-32}$ & $1.55\mathrm{e}{-33}$\\ 
 \hline
 $n=6$ & $2.42\mathrm{e}{-32}$ & $8.99\mathrm{e}{-34}$ & $6.00\mathrm{e}{-33}$ & $4.37\mathrm{e}{-34}$\\
 \hline
 $n=7$ & $1.31\mathrm{e}{-31}$ & $3.50\mathrm{e}{-33}$ & $1.78\mathrm{e}{-32}$ & $1.29\mathrm{e}{-33}$\\
 \hline
 $n=8$ & $7.38\mathrm{e}{-31}$ & $1.05\mathrm{e}{-32}$ & $5.37\mathrm{e}{-32}$ & $3.17\mathrm{e}{-33}$\\
 \hline
 \end{tabular}
\caption{Validation Error}
\label{tab:validation}
\end{table}

\section*{Potential Application: Moment-Derived Testing for p-value}

One potential application of these formulas, once placed in a proper statistical framework, is towards the development of a more accurate and computationally efficient estimation method of p-value of Pearson's correlation, as compared with existing approximation methods and permutation testing algorithms. 
\\
\\
Using such methods of estimation, one can obtain a valid approximation of the CDF with sufficiently less effort, especially since the distribution itself has no known tractable closed-form solution, while the moments are relatively inexpensive to obtain, depending on the moments of the data alone. Given a sufficient number of estimates of the moments, one can estimate the distribution of Pearson's correlation to any desired degree of accuracy. Such an analysis could potentially be more efficient, certainly more efficient than computing all the permutations of the correlation coefficient.  We posit that such a method exists for determining a direct p-value estimate.

\bibliographystyle{plain}
\bibliography{main.bib}


\end{document}